\documentclass[11pt,reqno]{amsproc}
\usepackage[all]{xy}
\usepackage[colorlinks,urlcolor=blue,citecolor=blue,linkcolor=blue]{hyperref}
\usepackage{mathpazo,natbib,xcolor,graphicx,bm,anysize}
\usepackage{caption,subcaption,eqnarray}

\pagecolor{lightgray}

\marginsize{15.3mm}{15.3mm}{6.3mm}{5.25mm}

\theoremstyle{theorem}
\newtheorem{theorem}{Theorem}[section]
\newtheorem{proposition}[theorem]{Proposition}
\newtheorem{lemma}[theorem]{Lemma}

\theoremstyle{remark}
\newtheorem*{pf}{\it Proof}
\newtheorem{remark}[theorem]{\it Remark}

\theoremstyle{definition}

\numberwithin{equation}{section}
\allowdisplaybreaks

\begin{document}
\title[Schr\"{o}dinger--Poisson--Slater equations subscaled near zero]{Schr\"{o}dinger--Poisson--Slater equations with nonlinearity subscaled near zero}
\author[S. Liu \& K. Perera]{Shibo Liu, Kanishka Perera\vspace{-1em}}
\dedicatory{Department of Mathematical Sciences, Florida Institute of Technology\\
Melbourne, FL 32901, USA}
\thanks{\emph{Emails}: \texttt{\bfseries lausb4@gmail.com} (S. Liu), \texttt{\bfseries kperera@fit.edu} (K. Perera)}
\begin{abstract}
We study the following zero-mass Schr{\"o}dinger--Poisson--Slater equation
  \[ - \Delta u + \left( \frac{1}{4 \pi | x |} \ast u^2 \right) u = f (| x |,
     u) \text{, \qquad} u \in \mathcal{D}^{1, 2} (\mathbb{R}^3) \text{} \]
  with nonlinearity subscaled near zero in the sense that $f (| x |, t)
  \approx a | t |^{p - 2} t$ as $| t | \rightarrow 0$ for some $p\in\big(\frac{18}{7},3\big)$. A
  nonzero solution is obtained via Morse theory when the nonlinearity is
  asymptotically scaled at infinity. For this purpose we prove an abstract
  result on the critical groups at infinity for functionals satisfying the geometric
  assumptions of the scaled saddle point theorem of Mercuri \& Perera [arXiv:2411.15887]. For
  the case that $f (| x |, \cdot)$ is odd, a sequence of solutions are
  obtained via a version of Clark's theorem due to Kajikiya [J.\ Funct.\ Anal.\ 225 (2005) 352--370].
\end{abstract}
\maketitle

\section{Introduction}

We consider the following zero-mass Schr{\"o}dinger--Poisson--Slater equation
\begin{equation}
  - \Delta u + \left( \frac{1}{4 \pi | x |} \ast u^2 \right) u = f (| x |, u)
  \text{, \qquad} u \in \mathcal{D}^{1, 2} (\mathbb{R}^3) . \label{e1}
\end{equation}
This equation arises from the study of nonlocal nonlinear Schr{\"o}dinger
equation
\begin{equation}
  i \partial_t \psi = - \Delta \psi + V (x) \psi + \left( \frac{1}{4 \pi | x
  |} \ast | \psi |^2 \right) \psi - \tilde{f} (| x |, | \psi |) \psi
  \text{\quad in } \mathbb{R}^3 \times \mathbb{R}_+ \label{SPs}
\end{equation}
in natural units, and its stationary counterpart. The equation (\ref{SPs})
comes from an approximation of the Hartree--Fock model of a quantum many-body
system of electrons, in which $| \psi |^2$ is the density of electrons and the
\ nonlocal convolution term represents the Coulombic repulsion between the
electrons. The local term $f (| x |, t) = | t |^{q - 2} t$ with $q =
\frac{8}{3}$ was introduced by  {\cite{PhysRev.81.385}} as a local
approximation of the exchange potential in the Hartree--Fock model (see, e.g.,
 {\cite{MR2013491}}, 
{\cite{MR1702877}}, and  {\cite{MR1836081}}).

The equation (\ref{e1}) has captured great interest in recent years. Partially
motivated by investigating the singular limit of
\[ - \Delta u + u + \lambda \left( \frac{1}{4 \pi | x |} \ast u^2 \right) u =
   | u |^{q - 2} u \text{\qquad in } \mathbb{R}^3 \]
as $\lambda \rightarrow 0$,  {\cite{MR2679375}} studied (\ref{e1}) with $f
(| x |, u) = | u |^{q - 2} u$ for $q \in (2, 3)$. The case that $q \in (3, 6)$
was studied by {\cite{MR2902293}}. These results were extended
by {\cite{Mercuri2016}}, where a more general equation
\[ - \Delta u + (I_{\alpha} \ast | u |^p) | u |^{p - 2} u = | u |^{q - 2} u
   \text{\qquad in } \mathbb{R}^N \]
being $I_{\alpha} (x) = A_{\alpha} | x |^{\alpha - N}$ the Riesz potential of
order $\alpha$, is considered. Based on the variational framework (the
Coulomb--Sobolev space $E  (\mathbb{R}^3)$ investigated in
{\cite{MR2679375}}) and convergence results (e.g., the nonlocal Brezis--Lieb lemma {\cite[Proposition 4.1]{Mercuri2016}}) established in these papers, many
subsequence results on (\ref{e1}) have been obtained.

{\cite{MR3912770}} studied the following equation with Sobolev
critical term $| u |^4 u$,
\[ - \Delta u + \left( \frac{1}{4 \pi | x |} \ast u^2 \right) u = \mu | u |^{p
   - 2} u + | u |^4 u \text{\qquad in } \mathbb{R}^3 \text{.} \]
By perturbation approach and truncation technique, they obtained a positive
ground state solution for $p \in (3, 4]$ with $\mu \gg 1$ or for $p \in (4,
6)$, and positive radial solution for $p \in \left( \frac{18}{7}, 3 \right)$
with $\mu \ll 1$, respectively. Their results for $p \in (3, 6)$ were improved
via simpler method by {\cite{Gu2024}}.

On the other hand, Schr{\"o}dinger--Poisson--Slater equation with
Coulomb-Sob\-o\-lev critical term $| u | u$
\begin{equation}
  - \Delta u + \left( \frac{1}{4 \pi | x |} \ast u^2 \right) u = | u | u + \mu
  | u |^{p - 2} u \text{\qquad in } \mathbb{R}^3 \text{} \label{o3}
\end{equation}
being $p \in (3, 6)$, was studied by  {\cite{Lei2023}}. By
minimizing the energy functional on the Nehari-Pohozaev manifold, a positive
ground state solution was obtained. Their result was generalized by   {\cite{Zhang2025}}, where in (\ref{o3}) the power nonlinearity $\mu | u
|^{p - 2} u$ is replaces by a more general $f (u)$.

Now we go back to our problem (\ref{e1}). Because the range of $p$ permitting
the embedding $E_r (\mathbb{R}^3) \hookrightarrow L^p (\mathbb{R}^3)$ is
$\left( \frac{18}{7}, 6 \right]$, to apply variational method it is natural to
require that the radial nonlinearity $f \in C (\mathbb{R}_+ \times
\mathbb{R})$ satisfies the growth condition:
\begin{itemize}
  \item[$(f)$] for some $p_{\pm} \in \left( \frac{18}{7}, 6 \right]$, there
  holds
  \[ | f (| x |, t) | \leq a_1 (| t |^{p_+ - 1} + | t |^{p_- - 1}) \text{.} \]
\end{itemize}
Then weak solutions of (\ref{e1}) are critical points of the $C^1$-functional
$\Phi : E_r \rightarrow \mathbb{R}$,
\[ \Phi (u) = \frac{1}{2} \int_{\mathbb{R}^3} | \nabla u |^2 \,\mathrm{d} x +
   \frac{1}{16 \pi} \iint_{\mathbb{R}^3 \times \mathbb{R}^3} \frac{u^2 (x) u^2
   (y)}{| x - y |}\, \mathrm{d} x \mathrm{d} y - \int_{\mathbb{R}^3} F (| x |, u) \,\mathrm{d}
   x \text{,} \]
where $F (| x |, t) = \int_0^t f (| x |, \cdot)$, $E_r = E_r (\mathbb{R}^3)$
is the radial Coulomb--Sobolev space given in (\ref{GS}). It is known that the sum $I (u)$ of the first two
terms in $\Phi$ has the scaling property $I (u_t) = t^3 I (u)$, where $u_t (x)
= t^2 u (t x)$ for $t > 0$. Recently,
{\cite{mercuri2024variationalmethodsscaledfunctionals}} exploited this scaling
feature and systematically developed a variational theory, including scaled
saddle point and linking theorems, and scaled local linking as well as related
critical group estimate. Applying to the problem (\ref{e1}), many new results
for various nonlinearity $f (| x |, t)$ were obtained in
{\cite{mercuri2024variationalmethodsscaledfunctionals}}.

The following scaled eigenvalue problem
\begin{equation}
  - \Delta u + \left( \frac{1}{4 \pi | x |} \ast u^2 \right) u = \lambda | u |
  u \text{, \qquad} u \in \mathcal{D}^{1, 2} (\mathbb{R}^3) \text{} \label{e2}
\end{equation}
is crucial to the study of (\ref{e1}). We say that $\lambda \in \mathbb{R}$ is
an eigenvalue of (\ref{e2}) if for this $\lambda$ the problem (\ref{e2}) has
nonzero solutions (called $\lambda$-eigenfunctions). It is clear that if $u$
is a $\lambda$-eigenfunction, so is the scaled function $u_t$, see
{\cite[\S5]{MR2902293}}. The set of all eigenvalues is denoted by $\sigma$. Using
$\mathbb{Z}_2$-cohomological index {\citep{MR0478189}}, a sequence of
eigenvalues of (\ref{e2}) carrying critical group information was constructed
by {\cite{mercuri2024variationalmethodsscaledfunctionals}}.

We start with the autonomous and asymptotically scaled case, that is $f (| x
|, t) = f (t)$ and
\begin{equation}
  \lim_{| t | \rightarrow \infty} \frac{f (t)}{| t | t} = \lambda \in
  \mathbb{R} \backslash \sigma \text{.} \label{If}
\end{equation}
Then, using the scaled saddle point theorem ({\citet[Theorem
1.5]{mercuri2024variationalmethodsscaledfunctionals}}), it has been shown in
{\cite[Theorem 1.18]{mercuri2024variationalmethodsscaledfunctionals}} that
(\ref{e1}) has a solution provided condition (\ref{E3}) below is satisfied.
With suitable subscaled conditions on $f$ near $t = 0$, we can obtain a nonzero
solution.

\begin{theorem}
  \label{t1}Let $f \in C (\mathbb{R})$, $\lambda \in \mathbb{R} \backslash
  \sigma$,
  \begin{equation}
    | f (t) - \lambda | t | t | \leq a_2 (| t |^{q_+ - 1} + | t |^{q_- - 1})
    \label{E3}
  \end{equation}
  for some $a_2 > 0$, $q_{\pm} \in \left( \frac{18}{7}, 3 \right)$. If
  \begin{equation}
    \lim_{| t | \rightarrow 0} \frac{F (t)}{| t |^3} = + \infty \text{,
    \qquad} 3 F (t) > f (t) t \geq 0 \text{\quad  for } t \neq 0 \text{,}
    \label{e3}
  \end{equation}
  then problem \eqref{e1} with $f (| x |, u) = f (u)$ has a nonzero solution.
\end{theorem}

As a simple example, we see that (\ref{e1}) with
\[
f (| x |, t) = \lambda | t |
t + | t |^{1 / 2} t
\]
has a nonzero solution provided $\lambda \in \mathbb{R}
\backslash \sigma$. If the limit in (\ref{e3}) holds we say that $f$ is
subscaled near zero. This is called $3$-sublinear near zero in the study of
$3$-Laplacian equations driven by $\Delta_3 u = \operatorname{div} (| \nabla u |
\nabla u)$.

It should be noted that the subscaled assumption (\ref{e3}) implies that $u=0$ is not a local minimizer of $\Phi$ and the mountain pass theorem is not applicable. The nonzero solution in Theorem \ref{t1} is obtained by comparing the critical
groups of $\Phi$ at $u = 0$ and at infinity. We prove an
abstract result on $C^* (\Phi, \infty)$, the critical groups at infinity, see Theorem \ref{t4}, which is a generalization of the classical result of
{\citet[Proposition 3.8]{MR1420790}}.

When $f (| x |, \cdot)$ is odd and subscaled near zero, we obtain a sequence of solutions without
any assumptions on $f$ except the subcritical growth condition $(f)$.

\begin{theorem}
  \label{t3}If $(f)$ is satisfied, $f (| x |, \cdot)$ is odd for all $x \in
  \mathbb{R}^3$ and
  \begin{equation}
    \lim_{| t | \rightarrow 0} \frac{F (| x |, t)}{| t |^3} = + \infty
    \text{\qquad uniformly in } x \in \mathbb{R}^3 \text{,} \label{Ff}
  \end{equation}
  then problem \eqref{e1} has a sequence of solutions $u_n$ such that $\Phi
  (u_n) \leq 0$ and $u_n \rightarrow 0$ in $E_r$.
\end{theorem}

\begin{remark}
  Theorem \ref{t3} improves {\cite[Theorem
  1.30]{mercuri2024variationalmethodsscaledfunctionals}}, which requires
  \[ f (| x |, t) = | t |^{\sigma - 2} t + g (| x |, t) \text{,} \]
  where $\sigma \in \left( \frac{18}{7}, 3 \right)$, the odd perturbation $g
  (| x |, \cdot)$ satisfies $(f)$ and $| G (| x |, t) | \leq C | t
  |^{\tilde{\sigma}}$ for some $\tilde{\sigma} \in \left( \sigma, 3
  \right)$.
  
 {\cite{Yang2019}} also obtained a sequence of negative energy
  solutions for the case that
  \begin{equation} f (x, t) = \mu k (x) | u |^{p - 2} u + | u |^4 u \text{,} \label{AQ}
  \end{equation}
  where $p \in (1, 2)$, $k \in L^{6 / (6 - p)} (\mathbb{R}^3)$ such that the
  Lebesgue measure $m (\{ k > 0 \}) > 0$, and $\mu \ll 1$.
\end{remark}

In our final result, we deal with the case that $f (| x |, t) = | t |^{q -
2} t - | t |^{p - 2} t$ with $q < \min \{ 3, p \}$. Without assuming $p \leq
6$ to ensure that $\Phi$ is well defined on $E_r$, we obtain a sequence of solutions.

\begin{theorem}
  \label{2t}If $q \in \left( \frac{18}{7}, 3 \right)$, $p \in (q, \infty)$,
  then the problem
  \begin{equation}
    - \Delta u + \left( \frac{1}{4 \pi | x |} \ast u^2 \right) u = | u |^{q -
    2} u - | u |^{p - 2} u \text{, \qquad} u \in \mathcal{D}^{1, 2} (\mathbb{R}^3)
    \label{es}
  \end{equation}
  has a sequence of solutions $u_n$ such that $\varphi (u_n) \leq 0$, $u_n
  \rightarrow 0$ in $E_r$, where $\varphi$ is the functional given in
  \eqref{f0}.
\end{theorem}

This result improves {\cite[Corollary
1.33]{mercuri2024variationalmethodsscaledfunctionals}}, where $p \in (q, 6]$
is required. Note that the case $f (| x |, t) = | t |^{2 / 3} t - | t |^{4 /
3} t$ arising in Thomas--Fermi--Dirac--von Weizs{\"a}cker model of Density
Functional Theory (see {\cite{MR629207,MR3762278}}) is covered in
both {\cite[Corollary 1.33]{mercuri2024variationalmethodsscaledfunctionals}}
and our Theorem \ref{2t}.

The paper is organized as follows. In Section \ref{SS2} we recall the radial
Coulomb--Sobolev space ${E_r}  (\mathbb{R}^3)$ and briefly discuss the
$s$-scaling on a Banach space $W$ (for us, $s = 3$ and $W = E_r
(\mathbb{R}^3)$), then prove our abstract result on $C^{*} (\Phi, \infty)$, the critical groups at infinity of $\Phi$.
Theorems \ref{t1}, \ref{t3} and \ref{2t} are proved in Sections \ref{SS3},
\ref{SS4} and \ref{SS5}, respectively.

\section{Preliminaries}\label{SS2}

Let $\mathcal{D}_r^{1, 2} (\mathbb{R}^3)$ be the Sobolev space of radial $\mathcal{D}^{1,
2}$-functions, following {\cite{MR2679375}}, on the radial Coulomb--Sobolev space
\begin{equation}
  E_r = E_r (\mathbb{R}^3) = \left\{ u \in \mathcal{D}_r^{1, 2} (\mathbb{R}^3) \left|
  \iint_{\mathbb{R}^3 \times \mathbb{R}^3} \right. \frac{u^2 (x) u^2 (y)}{| x
  - y |} \,\mathrm{d} x \mathrm{d} y < \infty \right\} \label{GS}
\end{equation}
we equip the norm
\begin{align} \| u \| = \left[ \int_{\mathbb{R}^3} | \nabla u |^2\, \mathrm{d} x + \left(
   \iint_{\mathbb{R}^3 \times \mathbb{R}^3} \frac{u^2 (x) u^2 (y)}{| x - y |}
  \, \mathrm{d} x \mathrm{d} y \right)^{1 / 2} \right]^{1 / 2} \text{.} \label{n}
  \end{align}
Then it has been proved in {\cite[Theorem 1.2]{MR2679375}} that $(E_r, \|
\cdot \|)$ is a uniformly convex Banach space which is embedded in $L^q (\mathbb{R}^3)$
continuously for $q \in \left( \frac{18}{7}, 6 \right]$ and compactly for $q
\in \left( \frac{18}{7}, 6 \right)$.

We consider the functional $\Phi : E_r \rightarrow \mathbb{R}$,
\begin{equation}
  \Phi (u) = \frac{1}{2} \int | \nabla u |^2  +
  \frac{1}{16 \pi} \iint \frac{u^2 (x) u^2
  (y)}{| x - y |} - \int F (| x |, u) \text{,} \label{F}
\end{equation}
here and in what follows, unless stated explicitly, all integrals are taken over $\mathbb{R}^3$, all double integrals are taken against $(x,y)\in\mathbb{R}^3\times\mathbb{R}^3$. When necessary, to be clear we write $\mathrm{d} x$ and $\mathrm{d} y$ at the end of the integrals (see e.g.\ (\ref{wQ})).
It is well known that under our assumption $(f)$, $\Phi$ is of class $C^1$,
critical points of $\Phi$ are weak solutions of problem (\ref{e1}).

The sum $I (u)$ of the first two terms in $\Phi$ enjoys certain scaling
property
\[ I (u_t) = t^3 I (u) \text{, \qquad where } u_t (x) = t^2 u (t x) \text{.}
\]
We now exploit this feature in a more general setting, see
{\cite[\S1.2]{mercuri2024variationalmethodsscaledfunctionals}} for
more details.

Let $W$ be a reflexive Banach space, $s \in (0, \infty)$. An $s$-scaling on
$W$ is a bounded continuous map $W \times [0, \infty) \rightarrow W$, $(u, t)
\mapsto u_t$ satisfying
\begin{itemize}
  \item  $u_0 = 0$, $u_1 = u$, $(u_t)_{t'} = u_{t t'}$ and $(\tau u)_t = \tau
  u_t$ for $t, t' \in [0, \infty)$, $\tau \in \mathbb{R}$ and $u \in W$;
  
  \item as $t \rightarrow \infty$, $\| u_t \| = O (t^s)$ uniformly in $u$ on
  bounded sets.
\end{itemize}
Let $I \in C^1 (W)$ be a coercive even functional satisfying the following
$s$-scaling property
\begin{equation}
  I (u_t) = t^s I (u) \text{\quad for } u \in W \text{, } t \geq 0 \text{;
  \qquad} I (u) > 0 \text{\quad for } u \in W \backslash \{ 0 \} \text{.}
  \label{se}
\end{equation}
Let $\mathcal{M}= I^{- 1} (1)$. For $u \in W \backslash \{ 0 \}$ there is a
unique
\[ \alpha (u) = I^{- 1 / s} (u) > 0 \]
such that $I (u_{\alpha (u)}) = 1$, that is $u_{\alpha (u)} \in \mathcal{M}$.
Since $\alpha : W \backslash \{ 0 \} \rightarrow \mathbb{R}$ is continuous,
the map
\begin{equation}
  \eta : [0, 1] \times (W \backslash \{ 0 \}) \rightarrow W \backslash \{ 0 \}
  \text{, \qquad} \eta (t, u) = u_{(1 - t) + t \alpha (u)} \label{h}
\end{equation}
is also continuous. Thus $\eta (1, \cdot)$ is a deformation from $W
\backslash \{ 0 \}$ to $\mathcal{M}$. Hence
\begin{equation}
  H^q (\mathcal{M}) \cong H^q (W \backslash \{ 0 \}) = 0 \text{\qquad for all
  } q \in \mathbb{N} \label{X}
\end{equation}
because $\dim W = \infty$, where $\mathbb{N}=\mathbb{Z} \cap [0, \infty)$,
$H^{\ast}$ stands for the cohomology with $\mathbb{Z}_2$ coefficients.

For $\Phi \in C^1 (W)$ we denote by $\mathcal{K}$ the set of all critical
points of $\Phi$. If $u$ is an isolated critical point of $\Phi$ with $\Phi
(u) = c$, $\ell \in \mathbb{N}$, we call
\[ C^{\ell} (\Phi, u) : = H^q (\Phi^c, \Phi^c \backslash \{ u \}) \]
the $\ell$-th critical groups of $\Phi$ at $u$, where $\Phi^c = \Phi^{- 1} (-
\infty, c]$.

If $\Phi$ satisfies the $({PS})$ condition and $c < \inf_{\mathcal{K}}
\Phi$ for some $c \in \mathbb{R}$, following {\citet[Definition
3.4]{MR1420790}}, we call
\[ C^{\ell} (\Phi, \infty) := H^{\ell} (W, \Phi^c) \]
the $\ell$-th critical groups of $\Phi$ at infinity. Note that by the
deformation lemma, all such $\Phi^c$ are homotopically equivalent, therefore
$H^{\ell} (W, \Phi^c)$ is independent of $c$ and $C^{\ell} (\Phi, \infty)$ is
well defined.

A very convenient way for getting nontrivial critical points is to compare the
critical groups of $\Phi$ at trivial critical point ($u = 0$) and at infinity.

\begin{proposition}[{\citet[Proposition 3.6]{MR1420790}}]
  \label{p9}If $\Phi \in C^1 (W)$ satisfies the $({PS})$ condition and
  $C^{\ell} (\Phi, 0) \neq C^{\ell} (\Phi, \infty)$ for some $\ell \in
  \mathbb{N}$, then $\Phi$ has a nonzero critical point.
\end{proposition}

We close this section by proving our abstract result on $C^{\ast} (\Phi,
\infty)$. For symmetric subset $A \subset W \backslash \{ 0 \}$, we denote by
$i (A)$ the $\mathbb{Z}_2$-cohomological index of $A$, which was introduced by {\cite{MR0478189}}, see also \cite[\S 2.5]{MR2640827} for a quick introduction.

Let $W$ be a Banach space with an $s$-scaling, $I \in C^1 (W)$ be an coercive
even functional satisfying the $s$-scaling property (\ref{se}), and
$\mathcal{M}= I^{- 1} (1)$.

\begin{theorem}
  \label{t4}Let $\Phi \in C^1 (W)$ satisfy the $({PS})$ condition, $A_0$
  and $B_0$ be disjoint nonempty closed symmetric subsets of $\mathcal{M}$
  such that
  \begin{equation}
    k = i (A_0) = i (\mathcal{M} \setminus B_0) < \infty \text{.} \label{i1}
  \end{equation}
  Let $B = \{ u_t \mid u \in B_0, t \geq 0 \}$. If $\inf_{B \cup \mathcal{K}}
  \Phi > - \infty$, and
  \begin{equation}
    \lim_{t \rightarrow + \infty} \Phi (u_t) = - \infty \text{\qquad uniformly
    on } A_0 \text{,} \label{e4}
  \end{equation}
  then $C^k (\Phi, \infty) \neq 0$.
\end{theorem}

\begin{remark}
  If the geometric assumptions of Theorem \ref{t4} hold, we say that $\Phi$
  has a scaled local linking in dimension $k$ at infinity. These geometric assumptions are closely related to that of the scaled saddle point theorem of \citet[Theorem 1.5]{mercuri2024variationalmethodsscaledfunctionals}. Theorem \ref{t4}
  reduces to {\cite[Proposition 3.8]{MR1420790}} if $W = W_1 \oplus W_2$ with
  $\dim W_1 = k$ is equipped with the $1$-scaling $u_t = t u$, $\mathcal{M}$
  is the unit sphere, $A_0 = W_1 \cap \mathcal{M}$ and $B_0 = W_2 \cap
  \mathcal{M}$.
\end{remark}

\begin{pf}
Consider the commutative diagram
\[
\xymatrix{H^k(\mathcal{M})&H^k(\mathcal{M},A_0)\ar[l]&H^{k-1}(A_0)\ar[l]_{\delta^*}&H^k(\mathcal{M})\ar[l]\\
H^k(\mathcal{M})\ar[u]&H^k(\mathcal{M},\mathcal{M}\backslash B_0)\ar[u]_{i_1^*}\ar[l]&H^{k-1}(\mathcal{M}\backslash B_0)\ar[u]_{i^*}\ar[l]&H^k(\mathcal{M})\ar[u]\ar[l]}
\]
where the rows are portions of exact cohomology sequences of the pairs
  $(\mathcal{M}, A_0)$ and $(\mathcal{M}, \mathcal{M} \backslash B_0)$ and the
  vertical maps are induced by the inclusion $i : (\mathcal{M}, A_0)
  \rightarrow (\mathcal{M}, \mathcal{M} \backslash B_0)$. From (\ref{X}) we
  have $H^k (\mathcal{M}) = 0$, thus $\delta^{\ast}$ is injective. Thanks to
  \citet[Proposition 2.14(iv)]{MR2640827}, condition (\ref{i1}) implies that $i^{\ast}$ is
  nontrivial. Therefore, $i_1^{\ast}$ is also nontrivial.
  
  Pick $a < \inf_{B \cup \mathcal{K}} \Phi$. Then $\Phi^a \subset W \backslash
  B$. By (\ref{e4}) there is $R \gg 1$ such that
  \[ \Phi (u_t) \leq a \text{\quad for } u \in A_R : = \{ u_t \mid u \in A_0,
     t \geq R \} \text{,} \]
that is, $A_R\subset\Phi^a$.

From $H^*(W,W\backslash\{0\})=0$ and the exact cohomology sequence
\[
\xymatrix{H^{k+1}(W,W\backslash\{0\})&\ar[l]H^k(W\backslash\{0\},A_R)&\ar[l]_{ \ \quad i_2^*}H^k(W,A_R)&\ar[l]H^k(W,W\backslash\{0\})}
\]
for the triple $(W,W\backslash\{0\},A_R)$, it follows that
\[
i_2^*:H^k(W,A_R)\to H^k(W\backslash\{0\},A_R)
\]
induced by the inclusion $(W\backslash\{0\},A_R)\to(W,A_R)$, is an isomorphism. Similarly,
\[
i_3^*:H^k(W,W\backslash B)\to H^k(W\backslash\{0\},W\backslash B)
\]
induced by the inclusion $(W\backslash\{0\},W\backslash B)\to(W,W\backslash B)$, is also an isomorphism.

Using the homotopies  $\eta$ given in (\ref{h}) and $h:(t,u)\mapsto u_{(1 - t) + t R\alpha (u)}$, it is easy to see that
  \[ 
  j_1:(\mathcal{M}, A_0) \rightarrow (W\backslash\{0\}, A_R) \text{, \qquad} j_2:(\mathcal{M},
     \mathcal{M} \backslash B_0) \rightarrow (W\backslash\{0\}, W \backslash B)  \]
  induce isomorphisms in cohomology, where $j_1:u\mapsto u_R$ and $j_2$ is the inclusion\footnote{\ Let $\pi:u\mapsto u_{\alpha(u)}$ be the projection to $\mathcal{M}$. Then
\begin{align*}
(j_1\circ\pi)=h(1,\cdot)\simeq h(0,\cdot)=1_{W\backslash\{0\}}\text{,\qquad}(\pi\circ j_1)=1_\mathcal{M}\text{,}\\
(j_2\circ\pi)=\eta(1,u)\simeq \eta(0,\cdot)=1_{W\backslash\{0\}}\text{,\qquad}(\pi\circ j_2)=1_\mathcal{M}\text{.}
\end{align*}}. We have the following  commutative diagram in which all maps except $j_1^*$ are induced by inclusion:
\[
\hspace{-1.5cm}\xymatrix{
&H^k(\mathcal{M}, A_0)&H^k(W\backslash\{0\},A_R)\ar[l]_{j_1^*}&H^k(W,A_R)\ar[l]_{\ \ i_2^*}
&H^k
(W,\Phi^a)\ar[l]_{j^*}\\
&H^k(\mathcal{M},
     \mathcal{M} \backslash B_0)\ar[u]^{i_1^*}&H^k(W\backslash\{0\},W\backslash B)\ar[u]\ar[l]_{j_2^*}&H^k(W,W\backslash B)\ar[l]_{\ \ \ i_3^*}\ar[ur]_{i_4^*}\ar[u]
}
\]
Because $i_1^{\ast}$ is nontrivial, chasing the diagram from $H^k(W,W\backslash B)$ to $H^k(\mathcal{M}, A_0)$ along the homomorphisms with labels, we conclude
  \[ C^k (\Phi, \infty) = H^k (W, \Phi^a) \neq 0 \text{.}\qedhere\]

Because $i_1^{\ast}$ is nontrivial, chasing the diagram from $H^k(W,W\backslash B)$ to $H^k(\mathcal{M}, A_0)$ along the homomorphisms with labels, we conclude
  \[ C^k (\Phi, \infty) = H^k (W, \Phi^a) \neq 0 \text{.}\qedhere\]
\end{pf}

\section{Proof of Theorem \ref{t1}}\label{SS3}

Now we consider the functional $\Phi : E_r \rightarrow \mathbb{R}$ given in
(\ref{F}). We write
\begin{equation}
  I (u) = \frac{1}{2} \int | \nabla u |^2 +
  \frac{1}{16 \pi} \iint \frac{u^2 (x) u^2
  (y)}{| x - y |} \label{11}
\end{equation}
for $u \in E_r$. Then $I \in C^1 (E_r)$ as well, and $I (u_t) = t^3 I (u)$ for
all $t \geq 0$.

The following lemma about $C^{\ell} (\Phi, 0)$ is motivated by 
{\cite{MR1634579}} and  {\cite{MR1982676}}, where semilinear and
quasilinear elliptic boundary value problems on a bounded domain $\Omega
\subset \mathbb{R}^N$ are considered, respectively.

\begin{lemma}
  If $f$ satisfies (f) and \eqref{e3}, then $C^{\ell} (\Phi, 0) \cong 0$ for
  all $\ell \in \mathbb{N}$.
\end{lemma}

\begin{pf}
  Let $u \in \Phi^{- 1} (0, \infty)$, we have
  \begin{equation}
    \int F (u_t(x))\, \mathrm{d} x = \frac{1}{t^3} \int F (t^2 u (y)) \,\mathrm{d} y \text{,
    \quad} \Phi (u_t) = t^3 I (u) - \frac{1}{t^3} \int F (t^2 u (y)) \,\mathrm{d} y
    \text{.} \label{wQ}
  \end{equation}
  By (\ref{e3}) and the Fatou's lemma,
  \begin{align*}
    \varliminf_{t \rightarrow 0} \frac{1}{t^3} \int F (u_t)\, \mathrm{d} x & = 
    \varliminf_{t \rightarrow 0} \int \frac{F (t^2 u (y))}{t^6}\, \mathrm{d} y\\
    & \geq  \int \varliminf_{t \rightarrow 0} \frac{F (t^2 u (y))}{t^6} \,\mathrm{d}
    y = + \infty \text{.}
  \end{align*}
  Thus (using the expression of $\Phi (u_t)$ in (\ref{wQ}))
  \[ \lim_{t \rightarrow 0} \frac{\Phi (u_t)}{t^3} = \lim_{t \rightarrow 0}
     \left( I (u) - \frac{1}{t^3} \int F (u_t) \,\mathrm{d} x \right) = - \infty
     \text{.} \]
  We conclude that $\Phi (u_t) < 0$ for $t \ll 1$. If
  \begin{equation}
    \Phi (u_t) = t^3 I (u) - \frac{1}{t^3} \int F (t^2 u (y)) \,\mathrm{d} y \geq 0
    \text{,}
  \end{equation}
  we have
  \[ I (u) \geq \frac{1}{t^{6}} \int F (t^2 u) \text{.} \]
  Using the inequality in (\ref{e3})  we get by differentiating (\ref{wQ})
  \begin{align*}
    \frac{\mathrm{d}}{\mathrm{d} t} \Phi (u_t) & =  3 t^2 I (u) + \frac{3}{t^4} \int
    F (t^2 u) - \frac{1}{t^3} \int f (t^2 u) (2 t u)\\
    & \geq  \frac{6}{t^4} \int F (t^2 u) - \frac{2}{t^2} \int f (t^2 u) u\\
    & =  \frac{2}{t^4} \int (3 F (t^2 u) - (t^2 u) f (t^2 u)) > 0 \text{.}
  \end{align*}
  We conclude that for $u \in \Phi^{- 1} (0, \infty)$, there is a unique $t
  (u) \in (0, 1]$ such that $\Phi (u_{t (u)}) = 0$. Moreover, applying the
  implicit function theorem to the map $(t, u) \mapsto \Phi (u_t)$, we see
  that the function $u \mapsto t (u)$ is continuous.
  
  Now, we define $\eta : [0, 1] \times E_r \rightarrow E_r$ by
  \[ \eta (s, u) = \left\{ \begin{array}{ll}
       u & \text{if } \Phi (u) \leq 0 \text{,}\\
       u_{(1 - s) + s t (u)} & \text{if } \Phi (u) > 0 \text{.}
     \end{array} \right. \]
  By the continuity of $(u, t) \mapsto u_t$ and $u \mapsto t (u)$, it is clear
  that $\eta$ is continuous. Thus $\eta (1, \cdot)$ is a deformation from
  $E_r$ to $\Phi_0 = \eta (1, E_r)$, and we conclude
  \begin{align*}
    C^{\ell} (\Phi, 0) & =  H^{\ell} (\Phi^0, \Phi^0 \backslash \{ 0 \}) =
    H^{\ell} (E_r, E_r \backslash \{ 0 \}) = 0
  \end{align*}
  for all $\ell \in \mathbb{N}$ because $\dim E_r = \infty$.\qed
\end{pf}

\subsubsection*{Proof of Theorem \ref{t1}}Since $C^{\ell} (\Phi, 0) \cong 0$ for
all $\ell \in \mathbb{N}$, if $\Phi$ satisfies the $({PS})$ condition and
$C^k (\Phi, \infty) \neq 0$ for some $k \in \mathbb{N}$, then by Proposition
\ref{p9}, $\Phi$ will have a nonzero critical point and Theorem \ref{t1} is
proved.

That $\Phi$ satisfies $({PS})$ and the geometric assumptions of our
Theorem \ref{t4} have been proved in
{\cite{mercuri2024variationalmethodsscaledfunctionals}}. Thus Theorem \ref{t4}
gives $C^k (\Phi, \infty) \neq 0$, and the proof of Theorem \ref{t1} is
complete.

For the reader's convenience, we outline the arguments in
{\cite{mercuri2024variationalmethodsscaledfunctionals}} below. We define
$C^1$-functionals $J$ and $G$ on $E_r$ via
\begin{equation}
  J (u) = \frac{1}{3} \int | u |^3 \text{, \qquad} G (u) = \int \left( F (u) -
  \frac{1}{3} | u |^3 \right) \text{.} \label{12}
\end{equation}
Let $A = I'$, $B = J'$ and $g = G'$. Then $A, B \in \mathcal{A}_s$ with $s =
3$ (see \cite[Definition
2.1]{mercuri2024variationalmethodsscaledfunctionals} for the meaning of $\mathcal{A}_s$), and satisfy $(H_6)$--$(H_{12})$ in {\cite[\S
2.1]{mercuri2024variationalmethodsscaledfunctionals}}. Moreover, condition
(\ref{E3}) implies that for all $v \in E_r$,
\begin{equation}
  \langle g (u_t), v_t \rangle = o (t^3) \| v \| \text{\qquad  as } t
  \rightarrow \infty \label{gG}
\end{equation}
uniformly in $u$ on bounded sets, see {\cite[Lemma
3.3(iii)]{mercuri2024variationalmethodsscaledfunctionals}}. Thus, {\cite[Lemma
2.26]{mercuri2024variationalmethodsscaledfunctionals}} can be applied to our
functional
\[ \Phi (u) = I (u) - \lambda J (u) - G (u) \text{,} \]
and we conclude that $\Phi$ satisfies the $({PS})$ condition.

If $\lambda > \lambda_1$, since $\lambda \notin \sigma$ we may assume that
$\lambda \in \left( \lambda_k {, \lambda_{k + 1}}  \right)$ for some $k\in\mathbb{N}$. Let $\mathcal{M}=
I^{- 1} (1)$, $\psi : \mathcal{M} \rightarrow \mathbb{R}$ be defined by $\psi
(u) = \dfrac{1}{J (u)}$. Then, $\lambda \in (\lambda_k, \lambda_{k + 1})$
implies
\[ i (\psi^{\lambda_k}) = i (\mathcal{M} \backslash \psi_{\lambda_{k + 1}}) =
   k \text{,} \]
see {\cite[Theorem 1.3(iii)]{mercuri2024variationalmethodsscaledfunctionals}}, where $\psi^{\lambda_k}=\psi^{-1}(-\infty,\lambda_k]$ and $\psi_{\lambda_{k+1}}=\psi^{-1}[\lambda_{k+1},\infty)$.

We take $A_0 = \psi^{\lambda_k}$, $B_0 = \psi_{\lambda_{k + 1}}$. 
Using (\ref{gG}), it has been shown in the proof of {\cite[Theorem
2.25]{mercuri2024variationalmethodsscaledfunctionals}} that
\[ \Phi (u_t) \leq - t^3 \left( \frac{\lambda}{\lambda_k} - 1 + o (1) \right)
   \rightarrow - \infty \text{\qquad as } t \rightarrow + \infty \]
uniformly in $u \in A_0 = \psi^{\lambda_k}$; and $\Phi$ is bounded from below
on
\[ B = \{ u_t \mid u \in B_0, t \geq 0 \} \text{,} \]
hence on $B \cup \mathcal{K}$ (we may assume that the critical set
$\mathcal{K}$ is finite, otherwise $\Phi$ already has infinitely many critical
points). The geometric assumptions of our Theorem \ref{t4} have been verified.

If $\lambda < \lambda_1$, as discussed in the proof of {\cite[Theorem
2.25]{mercuri2024variationalmethodsscaledfunctionals}}, $\Phi$ is bounded from
below and attains its global minimum at some $u \in E_r$. We have $C^0 (\Phi,
u) \neq 0$. Thus $u$ is a nonzero critical point because $C^{\ell} (\Phi, 0)
\cong 0$ for all $\ell \in \mathbb{N}$.

\section{Proof of Theorem \ref{t3}}\label{SS4}

To prove Theorem \ref{t3}, we need the following version of the Clark theorem.
We write $\gamma (A)$ for the genus of any symmetric subset $A \subset W
\backslash \{ 0 \}$.

\begin{proposition}[{\citet[Theorem 1]{MR2152503}}]
  \label{p1}Let $W$ be a Banach space and $\Phi \in C^1 (W, \mathbb{R})$ be an
  even coercive functional satisfying the $({PS})_c$ condition for $c
  \leq 0$ and $\Phi (0) = 0$. If for any $k \in \mathbb{N}$ there is a
  symmetric subset $A_k \subset W \backslash \{ 0 \}$ such that $\gamma (A_k)
  \geq k$ and
  \begin{equation}
    \sup_{A_k} \Phi < 0 \text{,} \label{y}
  \end{equation}
  then $\Phi$ has a sequence of critical points $u_k \neq 0$ such that $\Phi
  (u_k) \leq 0$, $u_k \rightarrow 0$ in $W$.
\end{proposition}

\begin{remark}
  In the statement of {\citet[Theorem 1]{MR2152503}} it is assumed that $\Phi$
  satisfies $({PS})$, namely $({PS})_c$ for all $c \in \mathbb{R}$.
  However, checking the proof of that theorem, it is clear that only
  $({PS})_c$ for $c \leq 0$ is necessary.
\end{remark}

Because our assumptions on $f$ is quite weak, it is impossible to show that
$\Phi$ satisfies the $({PS})$ condition. Therefore, we follow the idea
initiated by {\cite{MR4162412}} and improved by {\cite{MR4803162}} to overcome this difficulty.

Let $\phi : [0, \infty) \rightarrow [0, 1]$ be a decreasing
$C^{\infty}$-function such that $| \phi' (t) | \leq 2$,
\[ \phi (t) = 1 \text{\quad for } t \in [0, 1]  \text{, \qquad} \phi (t) = 0
   \text{\quad for } t \in [2, \infty) \text{.} \]
We consider the truncated functional $\Psi : E_r \rightarrow \mathbb{R}$,
\[ \Psi (u) = I (u) - \phi (I (u))  \int F (| x |, u) \text{.} \]
The derivative of $\Psi$ is given by
\begin{align}
  \langle \Psi' (u), v \rangle & = \langle I' (u), v \rangle - \phi (I (u)) 
  \int f (| x |, u) v - \left( \int F (| x |, u) \right) \phi' (I (u)) \langle
  I' (u), v \rangle \nonumber\\
  & = \left( 1 - \left( \int F (| x |, u) \right) \phi' (I (u)) \right) 
  \langle I' (u), v \rangle - \phi (I (u))  \int f (| x |, u) v  \label{e0}
\end{align}
for $u, v \in E_r$.

\begin{lemma}
  \label{l2}The functional $\Psi$ is coercive and satisfies the
  $({PS})_c$ condition for $c < 1$.
\end{lemma}

\begin{pf}
  It is clear that there is $R > 0$ such that if $\| u \| \geq R$ then $I (u)
  \geq 2$, so $\phi (I (u)) = 0$ and $\Psi (u) = I (u)$. Thus, $\Psi$ is
  coercive.
  
  Now let $c < 1$ and $\{ u_n \}$ be a $({PS})_c$ sequence of $\Psi$. We
  claim that for $n \gg 1$ there holds
  \begin{equation}
    1 - \left( \int F (| x |, u_n) \right) \phi' (I (u_n)) \geq 1 \text{.}
    \label{Xs}
  \end{equation}
  In fact, for $n \gg 1$,
  \begin{itemize}
    \item if $I (u_n) < 1$, then $\phi' (I (u_n)) = 0$ and (\ref{Xs}) holds;
    
    \item if $I (u_n) \geq 1$, then
    \[ - \phi (I (u_n)) \int F (| x |, u_n) = \Psi (u_n) - I (u_n) < 1 - I
       (u_n) \leq 0 \text{.} \]
    Hence $\int F (| x |, u_n) \geq 0$ and (\ref{Xs}) holds again because
    $\phi' (I (u_n)) \leq 0$.
  \end{itemize}
  Because $\Psi$ is coercive, $\{ u_n \}$ is bounded. We may assume that $u_n
  \rightharpoonup u$ in $E_r$. From the compactness of the embedding $E_r
  \hookrightarrow L^s (\mathbb{R}^3)$ for $s \in \left( \frac{18}{7}, 6
  \right)$ and the subcritical growth condition $(f)$, we have
  \[ \int f (| x |, u_n) (u_n - u) \rightarrow 0 \text{.} \]
  Consequently,
  \begin{align*}
    &  \hspace{-2em} \left( 1 - \left( \int F (| x |, u_n) \right) \phi' (I (u_n))
    \right) \langle I' (u_n), u_n - u \rangle\\
    & =  \langle \Psi' (u_n), u_n - u \rangle + \phi (I (u_n)) \int f (| x
    |, u_n) (u_n - u) \rightarrow 0 \text{.}
  \end{align*}
  Now from (\ref{Xs}) we conclude $\langle I' (u_n), u_n - u \rangle
  \rightarrow 0$. Using {\cite[Lemma
  3.2]{mercuri2024variationalmethodsscaledfunctionals}}, we deduce $u_n
  \rightarrow u$ in $E_r$.\qed
\end{pf}

\subsubsection*{Proof of Theorem \ref{t3}}Because $U = I^{- 1} (- \infty, 1)$ is
an open set containing $0$, and $\Psi (u) = \Phi (u)$ for $u \in U$, if $\Psi$
has a sequence of critical points $u_n$ such that $u_n \rightarrow 0$ in
$E_r$, then for large $n$, $u_n$ will be critical points of $\Phi$ and Theorem
\ref{t3} is proved.

Thanks to Lemma \ref{l2}, to get such a sequence of critical points for $\Psi$
via Proposition \ref{p1}, it suffices to construct for each $k\in\mathbb{N}$ a symmetric subset $A_k
\subset W \backslash \{ 0 \}$ satisfying $\gamma (A_k) \geq k$ and
\[ \sup_{A_k} \Psi < 0 \text{.} \]

Unlike in {\cite{MR4162412}} and {\cite{MR4803162}}, because $I (u) \leq c \| u \|^3$ for $\| u \| \ll 1$ is not true, the required
$A_k$ could not be taken as a small $k$-dimensional sphere unless the limit in
(\ref{Ff}) is strengthen to
\[ \lim_{| t | \rightarrow 0} \frac{F (| x |, t)}{| t |^2} = + \infty \text{,}
\]
so that we can use $I (u) \leq c \| u \|^2$ for $\| u \| \ll 1$. This is why \cite{Yang2019} needs $p \in (1, 2)$ in (\ref{AQ}), see \cite[Eq.\ (3.4)]{Yang2019}.

Therefore, to construct the set $A_k$ we should exploit the
$3$-scaling feature $I (u_t) = t^3 I (u)$ of $I$ more thoroughly. 
Consider the continuous map $h : E_r \backslash \{ 0 \} \rightarrow
\mathcal{M}$ defined via
\[ h (u) = u_{\alpha (u)} \text{, \qquad where } \alpha (u) = I^{- 1 / 3} (u)
   \text{.} \]
Given $k \in \mathbb{N}$, let $S^k$ be a $k$-dimensional unit sphere in $E_r$.
Because $h (S^k)$ is compact,
\begin{equation}
  \inf_{u \in h (S^k)} \int | u |^3 > 0 \text{, \qquad} \sup_{u \in h (S^k)}
  \int | u |^6 < \infty \text{.} \label{W}
\end{equation}
So we may choose $a \gg 1$ such that
\begin{equation}
  1 - a \inf_{u \in h (S^k)} \int | u |^3 < 0 \text{.} \label{Q}
\end{equation}
By $(f)$ and (\ref{Ff}), there is $b > 0$ such that
\begin{equation}
  F (| x |, t) \geq a | t |^3 - b | t |^6 \label{W1}
\end{equation}
for all $(x, t) \in \mathbb{R}^3 \times \mathbb{R}$. From (\ref{W}) and
(\ref{Q}), there is $\tau \in (0, 1)$ small enough such that
\begin{equation}
  \delta := 1 - a \inf_{u \in h (S^k)} \int | u |^3 + b \tau^6 \sup_{u
  \in h (S^k)} \int | u |^6 < 0 \text{.} \label{W2}
\end{equation}
Now we define a map $\pi : S^k \rightarrow E_r \backslash \{ 0 \}$ via
\begin{equation}
  \pi (u) = u_{\tau \cdot \alpha (u)} = (u_{\alpha (u)})_{\tau} \text{} =
  \bar{u}_{\tau} \text{,} \label{pi}
\end{equation}
where $\bar{u} = h (u) \in \mathcal{M}$, then set $A_k = \pi (S^k)$. Since
$\pi$ is an odd continuous map, we get
\[ \gamma (A_k) \geq \gamma (S^k) = k + 1 \geq k \text{.} \]
If $v \in A_k$, then $v = u_{\tau}$ for some $u \in h (S^k)$. Moreover, $\Psi
(v) = \Phi (v)$ because for $\tau \in (0, 1)$ we have $I (v) = \tau^3 I (u) <
1$. Now, using (\ref{W1}) and (\ref{W2}) we have
\begin{align*}
  \frac{\Psi (v)}{\tau^3} = \frac{\Phi (v)}{\tau^3} & =  \frac{\Phi
  (u_{\tau})}{\tau^3} = 1 - \frac{1}{\tau^6} \int F (| \tau^{- 1} x |, \tau^2
  u)\\
  & \leq  1 - \frac{1}{\tau^6} \int (a | \tau^2 u |^3 - b | \tau^2 u |^6)\\
  & =  1 - a \int | u |^3 + b \tau^6 \int | u |^6 \leq \delta < 0 \text{.}
\end{align*}
We deduce $\sup_{A_k} \Psi < 0$. This concludes the proof of Theorem \ref{t3}.

\begin{remark}
  \label{R}To the best of our knowledge, in all applications of the classical Clark's theorem and its
  variants, including Kajikiya's theorem (see Proposition \ref{p1}),
  {\cite[Theorem 1.1]{MR3400440}} and {\cite[Proposition 2.2]{MR871998}}, the
  symmetric subset $A_k \subset W \backslash \{ 0 \}$ can be taken as a
  $k$-dimensional small sphere. 
\end{remark}

\section{Proof of Theorem \ref{2t}}\label{SS5}

Formally, solutions of (\ref{es}) are critical points of
\begin{align}
  \varphi (u) & =  I (u) - \frac{1}{q} \int | u |^q + \frac{1}{p} \int | u
  |^p \nonumber\\
  & =  \frac{1}{2} \int| \nabla u |^2 + \frac{1}{16
  \pi} \iint \frac{u^2 (x) u^2 (y)}{| x - y
  |} - \frac{1}{q} \int | u |^q + \frac{1}{p} \int | u |^p
  \text{.} \label{f0} 
\end{align}
However, since $p$ may be greater than the critical exponent $6$, we may have
$\int | u |^p = \infty$ for some $u \in E_r$, the functional $\varphi$ is not
well defined on $E_r$. Motivated by {\cite{MR1971024}}, where the $p$-Laplacian
equation
\[ - \Delta_p u - \lambda g (x) | u |^{p - 2} u = k (x) | u |^{q - 2} u - h
   (x) | u |^{s - 2} u \text{, \qquad} u \in \mathcal{D}^{1, p} (\mathbb{R}^N) \]
is considered, we introduce a new space as our variational framework.

Let $X$ be the completion of the set of all radial $C_0^{\infty}
(\mathbb{R}^3)$ functions under the norm
\[ \| u \|_p = \| u \| + | u |_p \text{,} \]
where $\|\cdot\|$ is the norm defined in (\ref{n}). Then $(X, \| \cdot \|_p)$ is a Banach space continuously embedded into $E_r$.
It is standard to verify that $\varphi : X \rightarrow \mathbb{R}$ is well
defined and of class $C^1$, critical points of $\varphi$ are weak solutions of
(\ref{es}). Thus, to prove Theorem \ref{2t} we only need to find a sequence of
critical points for $\varphi : X \rightarrow \mathbb{R}$.

\begin{lemma}
  $\varphi : X \rightarrow \mathbb{R}$ is coercive.
\end{lemma}

\begin{pf}
  Let $\{ u_n \} \subset X$ be a sequence satisfying
  \[ \| u_n \|_p = \| u_n \| + | u_n |_p \rightarrow \infty \text{,} \]
  we need to show that $\varphi (u_n) \rightarrow + \infty$.
  
  If $\sup_n | u_n |_p < \infty$, then $\| u_n \| \rightarrow \infty$. Hence
  $I (u_n) \rightarrow + \infty$. Set $v_n = (u_n)_{\alpha (u_n)}$. Then $v_n
  \in \mathcal{M}$,
  \[ u_n = ((u_n)_{\alpha (u_n)})_{\alpha^{- 1} (u_n)} = (v_n)_{\alpha^{- 1}
     (u_n)} \text{, \qquad} \alpha^{- 1} (u_n) = I^{1 / 3} (u_n) \text{.} \]
  Therefore, since
  \[ \sup_{v \in \mathcal{M}} \int | v |^q < + \infty \text{} \]
  and $q < 3$, we deduce
  \begin{align*}
    \varphi  (u_n) & \geq  \varphi_1 (u_n) : = I (u_n) - \frac{1}{q} \int |
    u_n |^q\\
    & =  I (u_n) - \frac{1}{q} \int | (v_n)_{\alpha^{- 1} (u_n)} |^q\\
    & =  I (u_n) - \frac{[\alpha^{- 1} (u_n)]^{2 q - 3}}{q} \int | v_n |^q\\
    & \geq  I  (u_n) - \frac{I^{(2 q - 3) / 3} (u_n)}{q} \sup_{v \in
    \mathcal{M}} \int | v |^q \rightarrow + \infty \text{.}
  \end{align*}
  The above argument also shows that $\varphi_1 (u_n)$ is bounded from below.
  Therefore, if $\sup_n | u_n |_p = + \infty$, up to a subsequence we have $| u_n
  |_p \rightarrow \infty$, hence
  \[ \varphi (u_n) = \varphi_1 (u_n) + \frac{1}{p} \int | u |^p \geq \inf_n
     \varphi_1 (u_n) + \frac{1}{p} \int | u_n |^p \rightarrow + \infty
     \text{.} \tag*{\qed}\]
\end{pf}

\begin{lemma}
  $\varphi$ satisfies the $({PS})$ condition.
\end{lemma}

\begin{pf}
  Let $\{ u_n \} \subset X$ be a $({PS})$ sequence of $\varphi$. Then $\{
  u_n \}$ is bounded in $X$. Hence, $\{ u_n \}$ is bounded in both $E_r$ and
  $L^p (\mathbb{R}^3)$. Up to a subsequence we have
  \[ u_n \rightharpoonup u \text{\quad in } E_r \text{, \qquad} u_n
     \rightharpoonup u \text{\quad in } L^p (\mathbb{R}^3) \text{, \qquad} u_n
     \rightarrow u \text{\quad a.e. in } \mathbb{R}^3 \text{,} \]
  and $u_n \rightarrow u$ in $L^q (\mathbb{R}^3)$ thanks to the compactness of
  $E_r \hookrightarrow L^q (\mathbb{R}^3)$. We deduce
  \[ \int | u |^{p - 2} u (u_n - u) \rightarrow 0 \text{, \quad} \int | u |^{q
     - 2} u (u_n - u) \rightarrow 0 \text{, \quad} \int | u_n |^{q - 2} u_n
     (u_n - u) \rightarrow 0 \text{,} \]
  as well as $\langle I' (u), u_n - u \rangle \rightarrow 0$ because $I' (u)$
  is a continuous linear functional on $E_r$. Consequently,
  \[ \langle \varphi' (u), u_n - u \rangle = \langle I' (u), u_n - u \rangle +
     \int | u |^{p - 2} u (u_n - u) - \int | u |^{q - 2} u (u_n - u)
     \rightarrow 0 \text{.} \]
  Note that this does not follow from $u_n \rightharpoonup u$ in $X$ because
  we don't know whether $X$ is reflexive. Combining with $\varphi' (u_n)
  \rightarrow 0$, we get
  \begin{align}
    o (1) & =  \langle \varphi' (u_n) - \varphi' (u), u_n - u \rangle
    \nonumber\\
    & =  \langle I' (u_n) - I' (u), u_n - u \rangle + \int (| u_n |^{p - 2}
    u_n - | u |^{p - 2} u) (u_n - u) \nonumber\\
    &  \hspace{4em} - \int (| u_n |^{q - 2} u_n - | u |^{q - 2} u) (u_n - u) \nonumber\\
    & =  \langle I' (u_n) - I' (u), u_n - u \rangle + \int (| u_n |^{p - 2}
    u_n - | u |^{p - 2} u) (u_n - u) + o (1) \text{.} \label{EZ} 
  \end{align}
  Since $u_n \rightharpoonup u$ in $E_r$, using {\cite[Lemma 2.3]{MR2902293}}
  we have
  \begin{align*}
  \iint \frac{u_n^2 (x) u_n (y) u (y)}{| x - y |} &\rightarrow \iint
     \frac{u^2 (x) u^2 (y)}{| x - y |} \text{,} \\
\iint \frac{u^2 (x) u
     (y) u_n (y)}{| x - y |} &\rightarrow \iint \frac{u^2 (x) u^2 (y)}{| x - y
     |} \text{.} 
\end{align*}
  Therefore,
  \begin{align*}
    &   \hspace{-1em}\langle I' (u_n) - I' (u), u_n - u \rangle\\
    & =  \int \nabla u_n \cdot \nabla (u_n - u) + \iint \frac{u_n^2 (x)
    u_n (y) (u_n (y) - u (y))}{4 \pi | x - y |}\\
    &  \hspace{4em} - \int \nabla u \cdot \nabla (u_n - u) - \iint \frac{u^2 (x) u (y)
    (u_n (y) - u (y))}{4 \pi | x - y |}\\
    & =  \int | \nabla (u_n - u) |^2 + \frac{1}{4 \pi} \iint \frac{u_n^2 (x)
    u^2_n (y)}{| x - y |} + \frac{1}{4 \pi} \iint \frac{u^2 (x) u^2 (y)}{| x -
    y |}\\
    &  \hspace{4em} - \frac{1}{4 \pi} \iint \frac{u_n^2 (x) u_n (y) u (y)}{| x - y |} -
    \frac{1}{4 \pi} \iint \frac{u^2 (x) u (y) u_n (y)}{| x - y |}\\
    & =  \int | \nabla (u_n - u) |^2 + \frac{1}{4 \pi} \iint \frac{u_n^2 (x)
    u^2_n (y)}{| x - y |} - \frac{1}{4 \pi} \iint \frac{u^2 (x) u^2 (y)}{| x -
    y |} + o (1) \text{.}
  \end{align*}
  Since $u_n \rightarrow u$ a.e.\ in $\mathbb{R}^3$, this and the Fatou's lemma
  yield
  \begin{align}
    &  \hspace{-1em} \varliminf_{n \rightarrow \infty} \langle I' (u_n) - I' (u), u_n - u
    \rangle \nonumber\\
    & \geq  \varliminf_{n \rightarrow \infty} \frac{1}{4 \pi} \iint \frac{u_n^2
    (x) u^2_n (y)}{| x - y |} - \frac{1}{4 \pi} \iint \frac{u^2 (x) u^2 (y)}{|
    x - y |} \geq 0 \text{.} \label{ZS} 
  \end{align}
  Because (the integrand is nonnegative)
  \[ \int (| u_n |^{p - 2} u_n - | u |^{p - 2} u) (u_n - u) \geq 0 \text{,} \]
  From (\ref{ZS}) and (\ref{EZ}) we get
  \[ \varliminf_{n \rightarrow \infty} \int (| u_n |^{p - 2} u_n - | u |^{p - 2}
     u) (u_n - u) = 0 \text{, \quad } \varliminf_{n \rightarrow \infty} \langle
     I' (u_n) - I' (u), u_n - u \rangle = 0 \text{.} \]
  Up to a further subsequence, we have $\langle I' (u_n), u_n - u \rangle
  \rightarrow 0$, hence $\| u_n - u \| \rightarrow 0$ by {\cite[Lemma
  3.2]{mercuri2024variationalmethodsscaledfunctionals}}; and thanks to the
  H\"{o}lder inequality
  \begin{align*}
    0 & \leq  \left[ \left( \int | u_n |^p \right)^{(p - 1) / p} - \left(
    \int | u |^p \right)^{(p - 1) / p} \right] \left[ \left( \int | u_n |^p
    \right)^{1 / p} - \left( \int | u |^p \right)^{1 / p} \right]\\
    & \leq  \int (| u_n |^{p - 2} u_n - | u |^{p - 2} u) (u_n - u)
    \rightarrow 0 \text{,}
  \end{align*}
  which implies $\int | u_n |^p \rightarrow \int | u |^p$. Since $u_n
  \rightharpoonup u$ in $L^p (\mathbb{R}^3)$, we have $u_n \rightarrow u$ in
  $L^p (\mathbb{R}^3)$ and
  \begin{align*}
    \| u_n - u \|_p & =  \| u_n - u \| + | u_n - u |_p \rightarrow 0 \text{.}
  \end{align*}
  That is, $u_n \rightarrow u$ in $X$.\qed
\end{pf}

\subsubsection*{Proof of Theorem \ref{2t}}Given $k \in \mathbb{N}$, using idea as
in the proof of Theorem \ref{t3} we construct symmetric subset $A_k \subset X
\backslash \{ 0 \}$ satisfying $\gamma (A_k) \geq k$ and
\[ \sup_{A_k} \varphi < 0 \text{.} \]
Let $\mathfrak{M}= \{ u \in X \mid I (u) = 1 \}$ and consider the map $h : X
\backslash \{ 0 \} \rightarrow \mathfrak{M}$ defined via
\[ h (u) = u_{\alpha (u)} \text{, \qquad where } \alpha (u) = I^{- 1 / 3} (u)
   \text{.} \]
Let $S^k$ be a $k$-dimensional unit sphere in $X$. By compactness of $h (S^k)$
we have
\[ \inf_{u \in h (S^k)} \int | u |^q > 0 \text{, \qquad} \sup_{u \in h (S^k)}
   \int | u |^p < \infty \text{.} \]
Because $p > q$ and $q < 3$, there is $\tau > 0$ small such that
\[ \delta := 1 - \frac{\tau^{2 q - 6}}{q} \inf_{u \in h (S^k)} \int | u
   |^q + \frac{\tau^{2 p - 6}}{p} \sup_{u \in h (S^k)} \int | u |^p < 0
   \text{.} \]
Define $\pi : S^k \rightarrow X \backslash \{ 0 \}$ as in (\ref{pi}), that is
$\pi (u) = (h (u))_{\tau}$. Let $A_k = \pi (S^k)$, then $\gamma (A_k) \geq k$.

If $v \in A_k$, then $v = u_{\tau}$ for some $u \in h (S^k)$, and we have
\begin{align*}
  \frac{\varphi (v)}{\tau^3} & =  \frac{\varphi (u_{\tau})}{\tau^3} = 1 -
  \frac{\tau^{2 q - 6}}{q} \int | u |^q + \frac{\tau^{2 p - 6}}{p} \int | u
  |^p \leq \delta < 0 \text{.}
\end{align*}
It follows that $\sup_{A_k} \varphi < 0$.

Since $\varphi$ is coercive and satisfies $({PS})$, using Proposition
\ref{p1} we get a sequence of critical points $u_n$ such that $\varphi (u_n)
\leq 0$ and $u_n \rightarrow 0$ in $X$. Theorem \ref{2t} is proved.

\end{document}